# Decomposition of neuronal assembly activity via empirical de-Poissonization[*]

**Werner Ehm**

*Institute for Frontier Areas of Psychology and Mental Health*
*Wilhelmstraße 3a, 79098 Freiburg, Germany*
**e-mail:** `ehm@igpp.de`

**Benjamin Staude**

*Computational Neuroscience Group, RIKEN Brain Science Institute*
*2-1 Hirosawa, Wako City, Saitama 351-0198, Japan*
**e-mail:** `staude@brain.riken.jp`

**Stefan Rotter**

*Institute for Frontier Areas of Psychology and Mental Health*
*& Bernstein Center for Computational Neuroscience, Freiburg, Germany*
**e-mail:** `rotter@biologie.uni-freiburg.de`

**Abstract:** Consider a compound Poisson process with jump measure $\nu$ supported by finitely many positive integers. We propose a method for estimating $\nu$ from a single, equidistantly sampled trajectory and develop associated statistical procedures. The problem is motivated by the question whether nerve cells in the brain exhibit higher-order interactions in their firing patterns. According to the neuronal assembly hypothesis (Hebb [13]), synchronization of action potentials across neurons of different groups is considered a signature of assembly activity, but it was found notoriously difficult to demonstrate it in recordings of neuronal activity. Our approach based on a compound Poisson model allows to detect the presence of joint spike events of any order using only population spike count samples, thus bypassing both the "curse of dimensionality" and the need to isolate single-neuron spike trains in population signals.

## 1. Introduction

The cell assembly hypothesis (Hebb [13]) postulates dynamically interacting groups of neurons as building blocks of representation and processing of information in the neocortex. Synchronization of action potentials across neuronal groups was suggested as a potential signature for active assemblies (von der

---

[*]Research supported by IGPP Freiburg and BMBF grant 01GQ0420 to BCCN Freiburg.

Malsburg [19]). In this scenario, the assembly members would exhibit specific higher-order correlations. It was shown that single nerve cells are sensitive to higher-order correlations in their input, and that they can react to such input in a correlation-specific way (Abeles [1], Diesmann et al. [8], Kuhn et al. [16], [17]). The statistical procedures for the analysis of higher-order correlations from multiple parallel spike train recordings suggested so far (Perkel et al. [23], Martignon et al. [20], [21], Grün et al. [9], [10], Nakahara and Amari [22], Gütig et al. [11], [12]) suffer from the "curse of dimensionality": the estimation of the necessary parameters (of order $2^N$ for a population of $N$ neurons) poses serious problems, mainly due to insufficient sample sizes. As a consequence, most attempts to directly observe cell assemblies have resorted to pairwise interactions (Brown et al. [6]). However, when relying on second-order approaches one cannot draw direct conclusions about the presence of higher-order effects.

Here we present a novel procedure for the statistical analysis of higher-order patterns in massively parallel spike trains that circumvents the above-mentioned problem. It involves two elements: (i) We assume a model for neuronal populations with a simple descriptive parametrization of higher-order effects. Specifically, we employ correlated point processes, where correlations of any order are induced by "inserting" appropriate patterns of synchronous spikes. (ii) We base our inference on spike counts extracted from a single population (multi-unit) spike train, rather than from multiple parallel (single-unit) spike trains. This leads to a parsimoniously parametrized univariate estimation problem circumventing the curse of dimensionality, which greatly reduces the demands with respect to the size of empirical samples. This strategy also alleviates some of the difficulties arising when single-unit spike trains have to be extracted from multi-unit recordings; cf. Section 7.

The mathematical framework is as follows. Consider $N$ neurons labelled $1, 2, \ldots, N$, and let $\mathcal{A}$ denote the collection of all non-empty subsets of the set $\overline{N} = \{1, 2, \ldots, N\}$. Associated with the neurons is a collection of point processes [pp's] that describe their joint activity. For each subset $a \in \mathcal{A}$ there is a pp $X_a(t)$, $t \geq 0$,[1] the jump times of which indicate *simultaneous* firing of all neurons $i \in a$. Thus at the lowest level of single neurons there are $N$ pp's $X_i(t)$, $t \geq 0$ ($i \in \overline{N}$), where $X_i(t)$ denotes the number of "single-cell" spikes occurring in the $i$-th neuron during time interval $[0, t]$. At the level of neuron pairs there are $\binom{N}{2}$ pp's $X_{\{ij\}}(t)$ counting the number of "two-cell" spikes that occur simultaneously in neurons $i$ and $j$ up till time $t$. Next, there are $\binom{N}{3}$ pp's $X_{\{ijk\}}(t)$ counting the number of "triple-cell" spikes occurring simultaneously in neurons $i, j, k$ up till time $t$. And so on.

The $n$-cell processes with $n \geq 2$ account for possible $n$-fold interactions between the single neurons. These processes are not assumed to be observable; only their superposition is. Our goal is to disentangle this mixture, i.e., to extract information about 2-fold, 3-fold, etc. neuron interactions from the aggregate activity of the whole population.

---

[1] We represent pp's by their associated (cumulative) counting processes rather than as random measures.

Concerning the underlying $n$-cell processes we make the following assumptions.

(A1) All processes $X_a(t)$, $a \in \mathcal{A}$, are statistically independent.
(A2) For each $a \in \mathcal{A}$, $X_a(t)$ is a right-continuous Poisson (point) process [ppp], with constant jump rate $\mu(a)$, and such that $X_a(0) = 0$.

The sum activity of the population is described by the process

$$Z(t) = \sum_{a \in \mathcal{A}} |a| X_a(t), \ t \geq 0. \tag{1.1}$$

Quantitiy $Z(t)$ counts the number of spikes that occur anywhere in the population up till time $t$. Note that a (unit) jump of $X_a(t)$ at time $t = \tau$, say, contributes $|a|$ spikes ($|a|$ = cardinality of subset $a$), hence makes $Z(t)$ jump upwards by $|a|$ units at time $t = \tau$. (By (A1), (A2), the probability equals zero that there is any $t \geq 0$ that is a jump time for more than one of the processes $X_a$, $a \in \mathcal{A}$.)

From representation (1.1) and our assumptions it follows that $Z(t)$, $t \geq 0$, is a compound Poisson process with jump measure $\nu$ given by a weighted sum of Dirac point masses $\delta_n$,

$$\nu = \sum_{n=1}^{N} \nu_n \delta_n, \quad \text{where for every } n \geq 1, \quad \nu_n = \sum_{a \in \mathcal{A}, \ |a|=n} \mu(a) \tag{1.2}$$

is the jump rate of the process $Y_n := \sum_{a \in \mathcal{A}, |a|=n} X_a$, the superposition of all $n$-cell processes. In terms of the latter, one also has the slightly more condensed representation

$$Z(t) = \sum_n n Y_n(t), \ t \geq 0. \tag{1.3}$$

The Poisson process underlying all jumps, $\sum_n Y_n$, is called the carrier process. Its jump rate is $\sum_n \nu_n = \nu_+$.

The activity of individual neurons can be described by pp's $W_i$, $i \in \overline{N}$, where $W_i$ represents the spike train at neuron $i \in \overline{N}$. In terms of the $n$-cell processes $X_a$ these single-unit spike trains are given by $W_i = \sum_{a \ni i} X_a$. The collection $\{W_i, i \in \overline{N}\}$ of all single-unit pp's provides a multivariate representation of the population activity by means of correlated ppp's: each $W_i$ is a Poisson pp with intensity $\lambda_i = \sum_{a \ni i} \mu(a)$; the correlation (or "interaction") between $W_i$ and $W_j$ ($i \neq j$) results from those $n$-cell processes $X_a$ ($n \geq 2$) which simultaneously contribute to $W_i$ and $W_j$, i.e., correspond to firing patterns $a \in \mathcal{A}$ including both neurons $i, j$. Note, finally, that the superposition of all single-unit processes should yield the total activity of the population, which in fact it does:

$$\sum_i W_i = \sum_i \sum_{a \ni i} X_a = \sum_a \sum_{i \in a} X_a = \sum_a |a| X_a = Z. \tag{1.4}$$

It has to be pointed out that the subsequent developments exclusively refer to the process $Z$ and its assumed compound Poisson character. The precise relations between, and properties of, the single-unit processes $W_i$ and the $n$-cell processes $X_a$ are technically irrelevant as long as they are compatible with the compound Poisson model. They do matter, however, in regard to the interpretation of the model and the results derived on its basis. The problem to be dealt with can now be stated as follows.

**Problem:** Suppose one is given the values of process $Z(t)$ at $L$ equidistant time points $t_l = lh$, $1 \leq l \leq L$, with spacing $h > 0$. The task is to estimate the jump measure $\nu$, or at least, the first $M$ rates $\nu_n$ ($n \leq M$), from these data.

Suitable estimates of the $\nu_n$ could be used to answer specific questions in regard to neural networks, such as: Are pairwise correlations sufficient to describe the network activity? Do there exist correlations of order $> 10$, say? Etc.

Note that continuous observation of process $Z(t)$ would allow registration of every single jump, hence make the problem quite simple. In practice, one depends on discrete sampling, in which case it is unknown how many jumps have contributed to the increment across an interval (unless the increment is $\leq 1$). The information left after discretization is contained in the increments of $Z(t)$ across the bins $(t_{l-1}, t_l]$. These bin counts represent poissonized versions of the jump measure, which is why we call our goal de-Poissonization. De-Poissonization as a powerful analytic method has found numerous applications in, e.g., information theory, the analysis of algorithms, and combinatorics (Jacquet and Szpankowski [14], Borodin et al. [5]), but we are not aware of applications similar to those presented here. For a technically related, continuous variant of our estimation problem see Jongbloed et al. [15].

The paper is organized as follows. In Section 2 we introduce the rate estimates using Fourier inversion. Their basic asymptotic properties are derived in Section 3. A more detailed study of the conditions for the validity of the underlying approximations is deferred to Appendix A. The conncection between our approach and empirical de-Poissonization is explained in Section 4. Asymptotics-based statistical procedures are outlined in Section 5 and illustrated by Monte Carlo simulations in Section 6. The paper concludes with a discussion of some neurobiological and statistical aspects of our approach, addressing particularly the possible consequences of violations of the compound Poisson assumption and possible further developments.

## 2. Estimating the rates $\nu_n$ — a Fourier-based proposal

Given the bin width $h > 0$, let $Z_1, \ldots, Z_L$ denote the increments of process $Z(t)$ across the successive sampled time points, $Z_l = Z(lh) - Z((l-1)h)$. (To ease notation, we suppress the parameter $h$ wherever feasible.) The $Z_l$ are independent and identically distributed random variables with the common characteristic function

$$\gamma_h(\theta) = \exp\left[ h \int \left(e^{ix\theta} - 1\right) \nu(dx) \right] = \exp\left[ h \sum_{n=1}^{N} \nu_n(e^{in\theta} - 1) \right]. \quad (2.1)$$

Our (first) approach to the estimation of the $\nu_n$ depends on the simple observation that the exponent of (2.1) is essentially a Fourier series from which its coefficients, the rates $\nu_n$, can be recovered by Fourier inversion. With $\nu_+ = \sum_{n=1}^{N} \nu_n$ we have

$$\sum_{n=1}^{N} \nu_n e^{in\theta} = \nu_+ + h^{-1} \log \gamma_h(\theta), \quad (2.2)$$

so by inverting (2.2) we get back the coefficients of interest,

$$\begin{aligned} \nu_n &= \frac{1}{2\pi} \int_{-\pi}^{\pi} \left( \nu_+ + h^{-1} \log \gamma_h(\theta) \right) e^{-in\theta} d\theta \quad (2.3) \\ &= \frac{1}{2\pi} \int_{-\pi}^{\pi} h^{-1} \log \gamma_h(\theta) e^{-in\theta} d\theta \quad (n = 1, \ldots, M). \end{aligned}$$

Eq. (2.3) suggests to estimate the $\nu_n$ by substituting the unkown characteristic function $\gamma_h$ by the empirical characteristic function [ecf] of the observed increments,

$$\widehat{\gamma}_h(\theta) = L^{-1} \sum_{l=1}^{L} e^{i\theta Z_l},$$

so that the $\nu_n$ are estimated as follows,

$$\widehat{\nu}_n = \frac{1}{2\pi} \int_{-\pi}^{\pi} h^{-1} \log \widehat{\gamma}_h(\theta) e^{-in\theta} d\theta \quad (n = 1, \ldots, M). \quad (2.4)$$

Note that the number of neurons, $N$, does not enter (2.4), hence need not be known. The exponential form of $\gamma_h$ immediately yields the expression for the log characteristic function. Matters are more involved with the empirical version $\widehat{\gamma}_h$ whose (continuously defined) complex logarithm may traverse different branches. Although such happens only with small probability if $L$ is large (by the consistency of the ecf as $L \uparrow \infty$, uniformly in $\theta \in [-\pi, \pi]$; see e.g. Prakasa Rao [24, Ch. 8.3]), related difficulties may well occur in finite samples. Ways of fixing them will develop from an alternative view of our estimation problem outlined in Section 4. It has to be emphasized that the estimates can assume negative values. We do not attempt to avoid this drawback, and rather consider such a case as indicative of a small value of the corresponding rate parameter.

Linear functionals $\lambda = \sum_n c_n \nu_n$ of the $\nu_n$, with finitely many real constants $c_n$, are straightforwardly estimated by substituting $\nu_n$ by their estimates,

$$\widehat{\lambda} = \sum_n c_n \widehat{\nu}_n = \frac{1}{2\pi} \int_{-\pi}^{\pi} h^{-1} \log \widehat{\gamma}_h(\theta) C(\theta) d\theta, \quad (2.5)$$

where in the last expression $C(\theta) = \sum_k c_k e^{-ik\theta}$ denotes the Fourier series of the $c_k$ (which appears when interchanging summation and integration). Of particular interest in our context are functionals of the form $\lambda_m = \sum_{n=m}^{N} \nu_n$. Here $\lambda_m$ represents the total activity of all $n$-cell processes of orders $n \geq m$, and the question is whether $\lambda_m = 0$ for all $m$ greater than some $m_0$ (and if so, what is $m_0$).

Since the number $N$ of neurons generally is unknown, it would be desirable to replace $\lambda_m$ by the infinite tail sum $\rho_m = \sum_{n=m}^{\infty} \nu_n$, with the understanding that $\nu_n = 0$ for all but finitely many $n$. This is easily achieved on noting that the Fourier series associated with $\lambda_m$ is $\sum_{n=m}^{N} e^{-in\theta} = (e^{-im\theta} - e^{-i(N+1)\theta})/(1 - e^{-i\theta})$ and using the Riemann-Lebesgue lemma, according to which $\lim_{N \to \infty} \int_{-\pi}^{\pi} f(\theta) e^{-iN\theta} d\theta = 0$ for every integrable function $f \in L^1[-\pi, \pi]$. Thus

$$\rho_m = \lim_{N \to \infty} \frac{1}{2\pi} \int_{-\pi}^{\pi} h^{-1} \log \gamma_h(\theta) \frac{e^{-im\theta} - e^{-i(N+1)\theta}}{1 - e^{-i\theta}} d\theta$$

$$= \frac{1}{2\pi} \int_{-\pi}^{\pi} h^{-1} \log \gamma_h(\theta) R_m(\theta) d\theta \quad \text{where} \quad R_m(\theta) = \frac{e^{-im\theta}}{1 - e^{-i\theta}}; \quad (2.6)$$

and at least on any event where $(\log \widehat{\gamma}_h(\theta))/(1 - e^{-i\theta}) \in L^1[-\pi, \pi]$ it makes sense to estimate $\rho_m$ by the expression

$$\widehat{\rho}_m = \frac{1}{2\pi} \int_{-\pi}^{\pi} h^{-1} \log \widehat{\gamma}_h(\theta) R_m(\theta) d\theta. \quad (2.7)$$

The integrability condition holds in fact with probability tending to 1 as $L \uparrow \infty$. For the moment it suffices to note that the pole of $R_m$ at $\theta = 0$ is counterbalanced by the zero at $\theta = 0$ of the *log* (empirical) characteristic function.

## 3. Asymptotics

The basic properties of the proposed estimates that are required for, e.g., test construction can be derived applying the usual methods of asymptotic statistics (e.g., Prakasa Rao [24]). This is done here in a heuristic fashion at first. Precise conditions under which the approximations can be expected to be valid are provided in Appendix A.

By the law of large numbers, the ecf $\widehat{\gamma}_h(\theta)$ is close to its expected value, $\gamma_h(\theta)$, with high probability as $L$ becomes large. Let $\xi_h(\theta) = \frac{\widehat{\gamma}_h(\theta)}{\gamma_h(\theta)} - 1$. The Taylor expansion

$$h^{-1} \log \widehat{\gamma}_h(\theta) = h^{-1} \log \gamma_h(\theta) + h^{-1} \log(1 + \xi_h(\theta))$$

$$= \sum_{k=1}^{N} \nu_k(e^{i\theta k} - 1) + h^{-1} \left(\xi_h(\theta) - \frac{1}{2}\xi_h(\theta)^2 + \ldots\right) \quad (3.1)$$

inserted into (2.4), with only the first term kept, gives the stochastic approximation

$$\widehat{\nu}_n - \nu_n \approx \frac{1}{2\pi} \int_{-\pi}^{\pi} h^{-1} \xi_h(\theta) e^{-in\theta} d\theta =: U_n. \quad (3.2)$$

From (3.2) one readily obtains the approximate (normal) distribution of the estimates $\widehat{\nu}_n$ on calculating moments of $U_n$ and applying the central limit theorem. Clearly $\xi_h$ is Hermitean, i.e., $\overline{\xi_h(\theta)} = \xi_h(-\theta)$ for real $\theta$ (the bar denotes complex conjugation), and $E\xi_h(\theta) = 0$. Therefore $EU_n = 0$, and by the independence of the random variables $Z_l$ one has for real $\theta_1, \theta_2$

$$h^{-2} E\xi_h(\theta_1)\overline{\xi_h(\theta_2)} = \frac{1}{(hL)^2} \sum_{k,l} E\left(\frac{e^{i\theta_1 Z_k}}{\gamma_h(\theta_1)} - 1\right)\left(\frac{e^{-i\theta_2 Z_l}}{\gamma_h(-\theta_2)} - 1\right) \quad (3.3)$$

$$= \frac{1}{(hL)^2} \sum_l E\left(\frac{e^{i(\theta_1-\theta_2)Z_l}}{\gamma_h(\theta_1)\gamma_h(-\theta_2)} - 1\right)$$

$$= \frac{1}{h^2 L}\left(\frac{\gamma_h(\theta_1 - \theta_2)}{\gamma_h(\theta_1)\gamma_h(-\theta_2)} - 1\right)$$

$$= T^{-1}\Gamma_h(\theta_1, -\theta_2),$$

where $T = hL$ denotes "observation length," and $\Gamma_h$ the kernel

$$\Gamma_h(\theta_1, \theta_2) = \frac{1}{h}\left(\frac{\gamma_h(\theta_1 + \theta_2)}{\gamma_h(\theta_1)\gamma_h(\theta_2)} - 1\right). \quad (3.4)$$

Hence, using (3.2) along with the variable substitution $-\theta_2 \mapsto \theta_2$ in the integration, we find that the asymptotic covariance of estimates $\widehat{\nu}_m, \widehat{\nu}_n$ is given by

$$\text{ascov}(\widehat{\nu}_m, \widehat{\nu}_n) = EU_m\overline{U}_n = \frac{T^{-1}}{(2\pi)^2}\int_{-\pi}^{\pi}\int_{-\pi}^{\pi}\Gamma_h(\theta_1, \theta_2)\,e^{-im\theta_1 - in\theta_2}\,d\theta_1\,d\theta_2$$

$$\equiv T^{-1}\Omega_{m,n}. \quad (3.5)$$

Slightly more generally, the asymptotic covariance of the estimates $\widehat{\lambda}_b, \widehat{\lambda}_c$ of two linear functionals corresponding to finite vectors of constants $(b_k) = b, (c_k) = c$ is given by the bilinear form

$$T^{-1} b'\Omega c = T^{-1}\sum_{k,l}\Omega_{k,l}\,b_k\,c_l.$$

Alternatively, one may pass to the Fourier series

$$B(\theta) = \sum_k b_k e^{-ik\theta}, \quad C(\theta) = \sum_k c_k e^{-ik\theta}$$

and calculate the asymptotic covariance as

$$\text{ascov}(\widehat{\lambda}_b, \widehat{\lambda}_c) = \frac{T^{-1}}{(2\pi)^2}\int_{-\pi}^{\pi}\int_{-\pi}^{\pi}\Gamma_h(\theta_1, \theta_2)\,B(\theta_1)\,C(\theta_2)\,d\theta_1\,d\theta_2. \quad (3.6)$$

The latter procedure also applies to certain infinite compounds such as the tail sums $\rho_m$, in which case the Fourier series $B, C$ become $R_m, R_n$, respectively; see (2.7). (Note that $\Gamma_h(\theta_1, \theta_2)$ vanishes along the coordinate axes $\theta_1 = 0$ and $\theta_2 = 0$, cancelling the poles of $R_m(\theta_1)$ and $R_n(\theta_2)$.)

Some insight into the structure of the (co-)variances can be gained by letting the bin width $h$ tend to zero. Then by (3.4) and (2.1)

$$\Gamma_h(\theta_1, \theta_2) = \sum_{k=1}^{N}\nu_k\left(e^{i(\theta_1+\theta_2)k} - e^{i\theta_1 k} - e^{i\theta_2 k} + 1\right) + O(h).$$

When inserted into (3.5), this gives the approximation

$$\Omega_{m,n} = \nu_n\, \delta_{m,n} + O(h)$$

for ($T$ times) the asymptotic covariances of the $\widehat{\nu}_n$, which are thus uncorrelated to first order, with variance $\nu_n/T$. These approximative variances are optimal in the sense that they are identical with those applying to the case where each $n$-cell ppp is directly observable. This is in fact little surprising considering that if $h\nu_+$ is small, then any bin contains at most one jump of the carrier process with high probability. Note, finally, that $\Omega_{n,n} = \nu_n$ vanishes if $\nu_n = 0$, indicating that $\widehat{\nu}_n$ tends to zero at a faster rate than $O_p(T^{-1/2})$ in this case. Of course, these conclusions apply strictly only in the limit $h \downarrow 0$.

Let us end this section by indicating the kind of conditions that turn out, in Appendix A, to be decisive for the validity of the proposed normal approximations: *$T$ should be "large," $h$ and $\nu_+/T$ "small," and $h\nu_+$ is a critical quantity that should be "moderate" at most.*

## 4. Connection with analytic de-Poissonization

The estimates (2.4) can also be derived in a different manner using analytic de-Poissonization. This approach provides complementary information about our estimates that affords alternative ways of calculating the relevant quantities (other than by numerical integration), or even different estimation methods.

To achieve the necessary generality we will temporarily modify, and slightly abuse, the notation used so far. Let us first assume that

$$\phi(w) = E w^Z = \sum_{k=0}^{\infty} p_k w^k$$

is the generating function of some, for the moment arbitrary, discrete random variable $Z \geq 0$ with $P[Z = k] = p_k$ ($k \geq 0$). Suppose that $\phi$ is analytic in a neighborhood $W$ of the closed unit disc in the complex plane, and that it has no zeros within $W$. Then $\log \phi$, too, is analytic in $W$, and an application of Cauchy's integral formula[2] gives (on picking the main branch, where $\log 1 = 0$)

$$\frac{1}{2\pi i} \int_{|w|=1} \frac{\log \phi(w)}{w}\, dw = \log \phi(0) = \log p_0\,. \tag{4.1}$$

More generally we have for any $n \geq 0$

$$\frac{(\log \phi)^{(n)}(0)}{n!} = \frac{1}{2\pi i} \int_{|w|=1} \frac{\log \phi(w)}{w^{n+1}}\, dw = \frac{1}{2\pi} \int_{-\pi}^{\pi} \log \phi(e^{i\theta})\, e^{-i\theta n}\, d\theta\,. \tag{4.2}$$

Let $\beta_n$ denote the expression on the right-hand side of (4.2). Since $\phi^{(n)}(0)/n! = p_n$, (4.2) then yields the relations

$$\beta_0 = \log p_0,\ \ \beta_1 = \frac{p_1}{p_0},\ \ \beta_2 = \frac{p_2}{p_0} - \frac{1}{2}\left(\frac{p_1}{p_0}\right)^2,\ \ \beta_3 = \frac{p_3}{p_0} - \frac{p_2\, p_1}{p_0^2} + \frac{1}{3}\left(\frac{p_1}{p_0}\right)^3\ldots. \tag{4.3}$$

---

[2] For Cauchy's formula and other related facts see Ahlfors [2].

Returning to our setting we now let $Z$ stand for the generic specimen of the increments $Z_l$. Then $\phi(e^{i\theta}) = \gamma_h(\theta)$, and a comparison with (2.1) to (2.3) shows that the right-hand side of (4.2) becomes $h\nu_n$ for $n \geq 1$, and $-h\nu_+$ for $n = 0$. Together with the evaluation (4.3) this implies that *the rates $\nu_n$ are expressible as functions of the probabilities $p_k$, $0 \leq k \leq n$,*

$$h\nu_+ = -\log p_0, \quad h\nu_1 = \frac{p_1}{p_0}, \quad h\nu_2 = \frac{p_2}{p_0} - \frac{1}{2}\left(\frac{p_1}{p_0}\right)^2, \quad \ldots . \tag{4.4}$$

Moreover, since $\nu_+ = \rho_1$ we get for the tail sums $\rho_m$

$$h\rho_1 = -\log p_0, \quad h\rho_2 = h\rho_1 - h\nu_1 = -\log p_0 - \frac{p_1}{p_0}, \quad \ldots . \tag{4.5}$$

Taylor expansion of $\log \phi(w)$ about $w = 0$ is another way of arriving at these relations (Dayley and Vere-Jones [7, Sect. 5.2]); however, the use of Cauchy's formula fits better with the Fourier approach. The relations converse to (4.4),

$$p_0 = e^{-h\nu_+}, \quad p_1 = e^{-h\nu_+} h\nu_1, \quad p_2 = e^{-h\nu_+}\left((h\nu_1)^2/2 + h\nu_2\right), \quad \ldots \tag{4.6}$$

supply a basis for alternative estimation approaches to be explored elsewhere.

A key point of the preceding considerations is that the relations (4.1) to (4.5) do not only apply to the "true" quantities $\nu_n$, $p_k$, but likewise to the estimates $\widehat{\nu}_n$, the empirical frequencies $\widehat{p}_k = L^{-1}|\{l : Z_l = k\}|$, and the empirical characteristic function $\widehat{\gamma}_h(t) = \sum_{k \geq 0} \widehat{p}_k e^{itk}$ — provided the latter has winding number $\hat{n}_0 = 0$ w.r.t. the origin.[3] For in this case $\widehat{\gamma}_h(t) \equiv \widehat{\gamma}_h(e^{it})$ can be continued to an analytic function $\widehat{\gamma}_h(w) = \sum_k \widehat{p}_k w^k$ (a polynomial, in fact, since with probability 1 only finitely many $\widehat{p}_k$ are nonzero) that is zero-free within the unit disc, the latter by the argument principle. Thus, if the winding number is zero, the arguments that led to Eq. (4.2) are valid for the ecf, too, whence

$$\frac{(\log \widehat{\gamma}_h)^{(n)}(0)}{n!} = \frac{1}{2\pi}\int_{-\pi}^{\pi} \log \widehat{\gamma}_h(e^{i\theta}) e^{-i\theta n} d\theta = h\widehat{\nu}_n \qquad (n \geq 0). \tag{4.7}$$

The important conclusion (from the second equality) is that *the estimates $\widehat{\nu}_n$ defined by (2.4) and (4.7), respectively, are identical if $\hat{n}_0 = 0$.*

By evaluating the left-hand side of (4.7) one obtains as above

$$h\widehat{\nu}_+ = -\log \widehat{p}_0, \quad h\widehat{\nu}_1 = \widehat{p}_1/\widehat{p}_0, \quad h\widehat{\nu}_2 = \widehat{p}_2/\widehat{p}_0 - \frac{1}{2}(\widehat{p}_1/\widehat{p}_0)^2, \quad \ldots \tag{4.8}$$

if $\hat{n}_0 = 0$, showing that the rate estimates can be directly calculated from the histogram of the bin counts in this case. Explicit expressions for asymptotic (co-)variances can be derived in a similar way; cf. Appendix B. Unfortunately, formulae quickly become complicated, and a computer algebra system will be required to handle all but the most simple cases.

---

[3] The curve $[-\pi, \pi] \ni \theta \mapsto \widehat{\gamma}_h(\theta)$ traces out a closed loop within the unit disc of the complex plane. Loosely speaking, the winding number w.r.t. 0 counts how often that loop circles around 0 (counter-clockwise). For the precise definition see any text book about complex analysis, e.g. Ahlfors [2].

## 5. Tests

Tests of null-hypotheses pertaining to the parameters $\nu_k$ can be constructed on the basis of our estimates and (co-)variance approximations under the assumption that the estimates are consistent and approximately normally distributed. Practical implementation of any such test requires an estimate of the asymptotic (co-)variances, that is, of the kernel $\Gamma_h(\theta_1, \theta_2)$; the other quantities are known. One obvious possibility is to replace the "true" characteristic function $\gamma_h$ by its empirical version $\widehat{\gamma}_h$ everywhere in (3.4). For an alternative that is expected to be more stable statistically, note that by (3.4) and (2.1) the kernel $\Gamma_h(\theta_1, \theta_2)$ admits the explicit representation

$$\Gamma_h(\theta_1, \theta_2) = \frac{1}{h}\left(\exp\left[h \sum_{n=1}^{N} \nu_n \left(e^{i\theta_1 n} - 1\right)\left(e^{i\theta_2 n} - 1\right)\right] - 1\right), \qquad (5.1)$$

which allows us to estimate $\Gamma_h(\theta_1, \theta_2)$ by substituting estimates $\widehat{\nu}_n$ for the unknowns $\nu_n$. The $\widehat{\nu}_n$ can assume negative values, with potentially unpleasant consequences. To avoid such, we plug in positive parts $\widehat{\nu}_n^+$ instead ($x^+ = x$ if $x > 0$, and $x^+ = 0$ otherwise). In practice the sum in (5.1) has to be truncated, and the proposed estimate of $\Gamma_h(\theta_1, \theta_2)$ becomes

$$\widehat{\Gamma}_h(\theta_1, \theta_2) = \frac{1}{h}\left(\exp\left[h \sum_{n=1}^{K} \widehat{\nu}_n^+ \left(e^{i\theta_1 n} - 1\right)\left(e^{i\theta_2 n} - 1\right)\right] - 1\right) \qquad (5.2)$$

with some suitably chosen number $K$. For example, if a null-hypothesis implies $\nu_k = 0$ for $k > m$, one may take $K = m$. This choice is not mandatory, however. One also could set $K = M$ if $\nu_1, \ldots, \nu_M$ are the rates that are being estimated.[4] Whatever the concrete choice of $\widehat{\Gamma}_h$, let $T^{-1}\widehat{\Omega}_{m,n}$, $T^{-1}\widehat{\Sigma}_{m,n}$, ... denote the estimated covariances obtained by replacing $\Gamma_h$ in (3.5), (B.1), ... by $\widehat{\Gamma}_h$.

Standard ANOVA type tests impose linear restrictions on the rates such as $\mathcal{H}_0$: "$A\nu = 0$," where $\nu$ denotes the column vector with entries $\nu_1, \ldots, \nu_M$ and $A$ is some $q \times M$ matrix of full rank $q$. The corresponding Wald type test statistic $W = T(A\widehat{\nu})'(A\widehat{\Omega}A')^{-1}A\widehat{\nu}$ is approximately $\chi^2$-distributed with $q$ degrees of freedom under $\mathcal{H}_0$. Related two-sample analogs may be used to test for equality of $\nu$ in different subjects, or under different conditions in the same subject.

Other tests relevant to our goals pertain to the tail sums $\rho_m$. Let $\overline{\nu}$ denote the maximal index $k$ such that $\nu_k > 0$. Clearly then, $\rho_m > 0$ for every $m \leq \overline{\nu}$ and $\rho_m = 0$ for every $m > \overline{\nu}$. This suggests to consider null-hypotheses of the form $\mathcal{H}_0$: "$\rho_m = 0 \, \forall \, m_1 \leq m \leq m_2$," for some fixed pair $m_1, m_2$ such that $2 \leq m_1 \leq m_2$. An associated test would reject $\mathcal{H}_0$ if $V = \max_{m_1 \leq m \leq m_2} V_m$ exceeds some critical value, where $V_m = (T/\widehat{\Sigma}_{m,m})^{1/2}\widehat{\rho}_m$. An approximate null-distribution for $V$ is difficult to derive theoretically (except for the trivial case $m_1 = m_2$), but may be estimable by means of the bootstrap method. Alternatively, one may test each single hypothesis "$\rho_m = 0$" observing that $V_m$ is approximately standard normally distributed (if $\rho_m = 0$), and finally apply some test level correction. Most useful in applications may be a simple screening procedure: plotting the standardized statistics $V_m$ (or the corresponding p-values) against $m$ could suffice to reveal the major features of interest.

---

[4] Ultimately, the choice of $K$ should depend on considerations related to the thereby achieved validity and efficiency of the resulting test.

## 6. Simulation results

In the following we present simulation results for several parameter configurations, including one where the Fourier inversion method fails to perform reliably. For Example 1 we chose $(\nu_1, \ldots, \nu_5) = (40, 10, 4, 3, 1)$ and $\nu_k = 0$ for all $k > 5$. Observation length was $T = 30$, and bin width $h = 0.02$. (All times are specified in seconds, all rates are given in events/second.) For Example 2 we put $\nu_1 = 150$ and $\nu_7 = 7$, and $\nu_n = 0$ otherwise. In this case we had $T = 60$ and $h = 0.005$. For each example we generated 50 realizations of the respective compound Poisson process by Monte Carlo simulation.

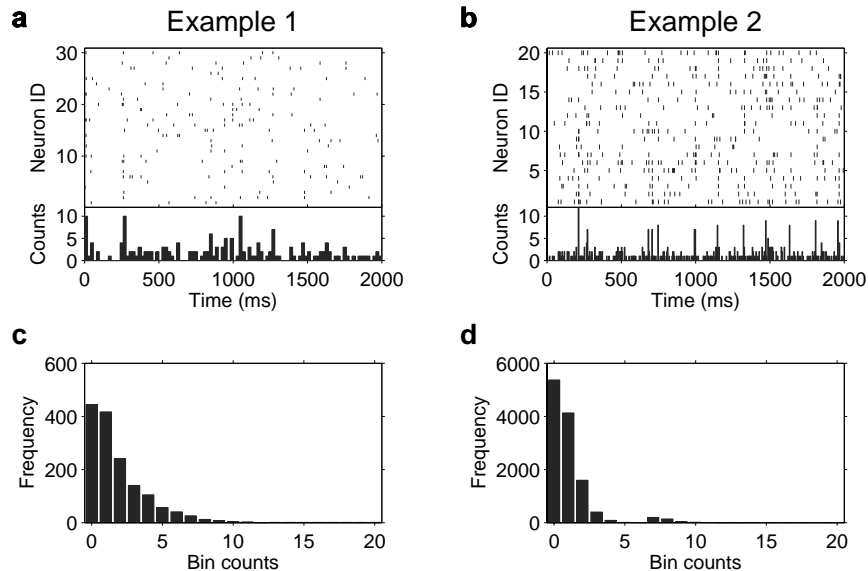

FIG 1. *Simulated examples.*

Graphical presentation of a numerical simulation of Example 1 (*a,c*) and Example 2 (*b,d*), see text for parameters. *a,b*, Raster plots, horizontal ticks mark the occurrence of a spike in the indicated neuron, some spikes are aligned into spike patterns (top). Discrete samples of the compound counting process (bin counts, bottom) shown for 2 seconds simulation time. *c,d*, Histogram of bin counts for a full-length simulation.

In order to illustrate the connection with the usual multivariate approach we split each process into a multivariate Poisson process representing a hypothetical population of neurons (cf. the paragraph preceding Eq. (1.4)). The populations comprised $N = 30$ and $N = 20$ neurons, respectively, in the two examples shown, yielding average single-neuron firing rates of about 3 Hz and 10 Hz. For simplicity, we chose a maximum entropy representation where all neurons have the same likelihood to participate in a particular pattern. That is, we generated independent ppp's $Y_n$ with intensities $\nu_n$ and assigned the $n$ spikes occurring at each jump of $Y_n$ randomly to one of the $\binom{N}{n}$ $n$-element subsets $a$ of the neurons. Recall, however, that the single spike trains merely are shown for illustration: our statistical analyses only depend on the bin counts of the compound process. These are displayed below the single spike trains in Figure 1.

For each simulated realization we estimated $\widehat{\nu}_n$ and $\widehat{\rho}_n$ for $n = 1, \ldots, 12$ using the Fourier method, with associated covariance matrix estimates obtained according to (5.2). As an illustration of hypothesis testing we calculated test power profiles $\beta_n$ as follows: for each $n = 1, \ldots, 12$ we determined $\beta_n$ as the relative frequency among the 50 simulations in which the test statistic $V_n = (T/\widehat{\Sigma}_{n,n})^{1/2} \widehat{\rho}_n$ exceeded the value 2, roughly corresponding to a one-sided 5 % level test of the null-hypothesis "$\rho_n = 0$". The estimated rate profiles reproduce the "true" profiles fairly well (Figure 2, left-hand column). The sharp decrease of the power profile in Figure 2d) properly reflects the vanishing of all $\nu_n$ for $n > 7$. Such a clear separation cannot be expected in Example 1 where the transition from the non-zero to the zero rates is rather smooth; see Figure 2c).

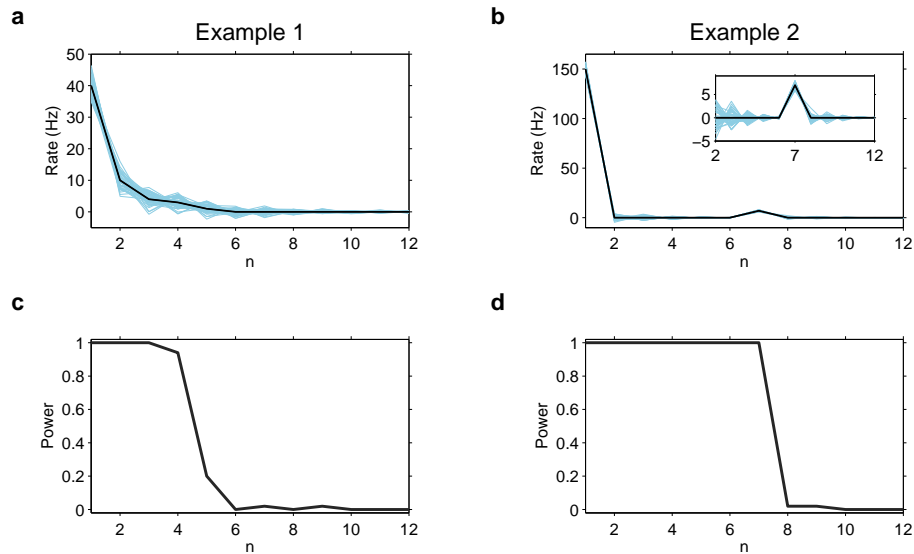

FIG 2. *De-Poissonization of sample data.*

Fourier-based de-Poissonization applied to the examples described in the text. *a,b,* Shown are the true ($\nu_n$, black) and estimated ($\widehat{\nu}_n$, blue) rates, plotted against jump height $n$ ("Amplitude"), for each of the 50 realizations. *c,d,* Test power profiles $\beta_n$ obtained by Monte Carlo simulation, corresponding to null-hypotheses "$\rho_n = 0$".

Simulations sometimes exhibit a peculiar phenomenon: the $\widehat{\nu}_n$ tend to form clusters. While in benign cases there is only one cluster concentrated at the true rates, increase of $h\nu_+$ induces splits of the estimated rates into two or more distinct clusters located at totally wrong rate profiles. The reason for this phenomenon is that the ecf has a nonzero winding number $\hat{n}_0$ in such a case, making the log-ecf end up in a different branch. This gives us a simple diagnostics: the phase angle of $\widehat{\gamma}_h(\theta)$, $\Im\log\widehat{\gamma}_h(\theta)$, then changes by a nonzero amount $2\pi\hat{n}_0$ along the loop $[-\pi,\pi] \ni \theta \mapsto \widehat{\gamma}_h(\theta)$. In particular, $\Im\log\widehat{\gamma}_h(\pi) - \Im\log\widehat{\gamma}_h(0)$ equals a nonzero multiple of $\pi$ in such a case (and is zero otherwise).

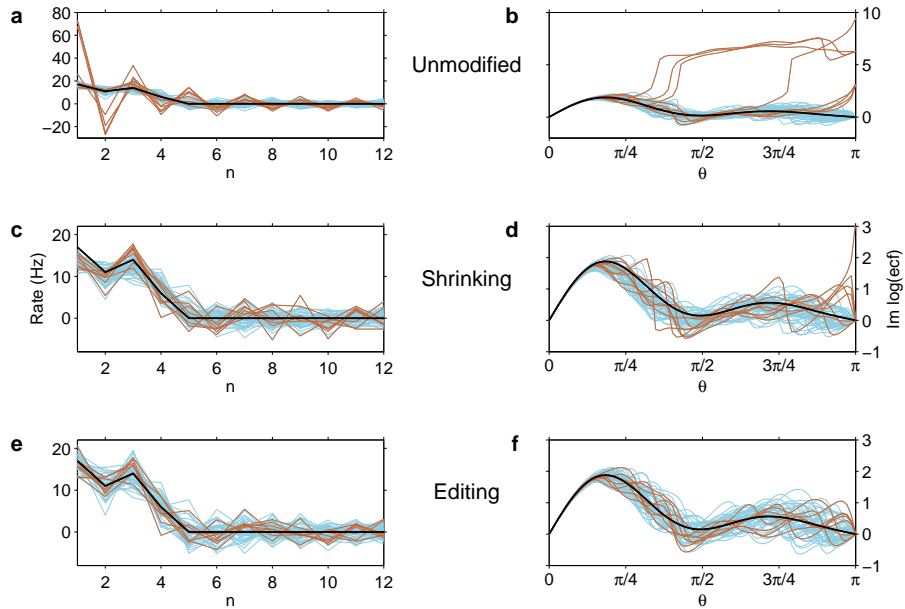

FIG 3. *Winding number problems solved by ecf-modification.*
An example where plain Fourier-based de-Poissonization fails. Rate profiles (*a,c,e*) and imaginary parts of the log characteristic functions (*b,d,f*) of the "true" model are represented by black lines, the empirical analogs by blue or red lines, respectively, if the winding number of the initial ecf is zero or non-zero. The deviating rate profiles (*a*) associated with a non-zero winding number of the ecf (*b*) are effectively corrected by "shrinking" (*c,d*) and "editing" (*e,f*), with one exception in case of shrinking.

A simple remedy is to shrink the loop towards 1 so that it does not circle around 0 anymore: given some $\delta \in (0,1)$, replace $\widehat{\gamma}_h(\theta)$ by $\widehat{\gamma}_h(\theta;\delta) = \delta + (1-\delta)\widehat{\gamma}_h(\theta)$, and obtain rate estimates $\widehat{\nu}_k(\delta)$ based on $\widehat{\gamma}_h(\cdot;\delta)$ (rather than on $\widehat{\gamma}_h(\cdot)$). As long as $\delta$ is small, the resulting perturbation of the estimates should be small, too. For a less crude corrective one starts out directly from the origin of non-zero winding numbers: such result if $\widehat{\gamma}_h$ considered as a polynomial in $w \in \mathbf{C}$ has zeros within the unit disc (cf. Section 4). The idea here is to move any such zero "outside" by multiplying it with a suitable factor. For definiteness, let $\epsilon > 0$ and suppose that $\widehat{\gamma}_h(w) = \sum_{k=0}^{K} \widehat{p}_k w^k$ ($K = \max_l Z_l$) has zeros $\alpha_1,\ldots,\alpha_K$, $\widehat{\gamma}_h(w) = \prod_{k=1}^{K}(w-\alpha_k)/(1-\alpha_k)$. Let "edited zeros" $\widetilde{\alpha}_k = \widetilde{\alpha}_k(\epsilon)$ be defined by $\widetilde{\alpha}_k = (1+\epsilon)\alpha_k/|\alpha_k|$ if $|\alpha_k| \leq 1+\epsilon$, and $\widetilde{\alpha}_k = \alpha_k$ otherwise[5], and put

$$\widetilde{\gamma}_h(w;\epsilon) = \prod_{k=1}^{K} \frac{w - \widetilde{\alpha}_k}{1 - \widetilde{\alpha}_k}.$$

Then define rate estimates $\widetilde{\nu}_n(\epsilon)$ as above, now replacing $\widehat{\gamma}_h(\cdot)$ by $\widetilde{\gamma}_h(\cdot;\epsilon)$.

Figure 3 illustrates both procedures for the following choice of parameters: $(\nu_1, \ldots, \nu_4) = (17, 11, 14, 6)$ and $\nu_k = 0$ for $k > 4$; $T = 60$; $h = 0.05$;[6] shrinking parameter $\delta = 0.02$; editing parameter $\epsilon = 0.075$. Our simulations indicate that editing induces less bias in the rate estimates than shrinking. Moreover, it works always, while for shrinking this holds true only if $\delta$ is selected adaptively. On the other hand, such an adaptive choice is easily implemented, and shrinking does not require the computation of all zeros of $\widehat{\gamma}_h$. We will come back to these matters elsewhere, along with more extensive simulation studies and applications to measured data.

## 7. Discussion

In this paper we have proposed a compound Poisson model along with associated statistical procedures for the analysis of higher-order interactions in neuronal assembly activity. The proposed model is purely phenomenological; explanation of neural mechanisms is neither intended nor within its scope. Its main purpose is to provide a formal framework for the statistical analysis of active assemblies in sampled data, to the degree they are reflected by the occurrence of higher-order spike patterns. This goal is achieved without recourse to simultaneous recordings of the single-neuron spike trains. Only the compound signal is used, which comprises the superimposed spike trains of all neurons. It is then impossible to identify specific groups of neurons that participate in specific coordinated events. In return, our approach allows to design statistical procedures that operate reliably on much smaller samples than other methods suggested previously. Moreover, it avoids certain technical problems associated with the identification of individual neurons in population signals typically obtained with extracellular recording electrodes: Although spikes have still to be isolated from noise, and simultaneous spikes on a single channel have to be correctly identified as such, our approach practically eliminates the need of "spike sorting", a difficult and error-prone procedure (often based on a classification of spike wave-forms) through which spikes are assigned to individual neurons.

Technically, our approach hinges on the compound Poisson character of the bin count distribution. Ultimately, it is this assumption that establishes the connection between the observed counts and the not directly observable higher-order correlations in the activity of the neuron population. Doing without any such functional relation is impossible in inferences based on the compound activity only. On the other hand, as with any model the assumptions may not be adequate in applications, and potentially serious consequences have to be taken into account if they are violated. Let us adress some issues that might turn out problematic with actual neuronal behavior.

---

[5]It may be useful to "edit" zeros outside but close to the unit disc as well. Other "zero editing" such as reflection at the unit circle may also be worth considering.

[6]The critical parameter $h\nu_+$ is 2.4 in this case, whereas it was about 1 in the two well-behaved examples discussed above.

It is currently hard to judge how well a real neuronal population conforms to the Poisson model, although irregular (i.e. Poissonian) spike trains are considered to be characteristic of normal operation of neocortical networks. Recently, evidence is accumulating that the single-neuron firing rates in the cortex of behaving animals are extremely low ($\ll 1$ Hz, Lee et al. [18]), much lower than previously thought (Shadlen and Newsome [26]). In our model, individual spikes are distributed over a large number of different assemblies, so "statistics of rare events" may be an adequate approximation. Deviations from Poisson statistics in real neurons, on the other hand, seem to occur mainly for higher firing rates, e.g. during stimulus response (Amarashingham et al. [3]) or in the large projection neurons of the motor cortex (Baker and Lemon [4]).

Our assumption of identically distributed bin counts means that the compound process $Z$, and the carrier process $Y_+$, respectively, should have stationary increments within the observation period $[0, T]$. While quasi-stationary stretches of data of sufficient length could possibly be obtained in anesthetized animals, the task-dependent modulations of neuronal firing rates in behaving animals generate inevitable non-stationarities. The consequences of unattended firing rate modulations are potentially quite serious: fluctuating rates increase the weight in the tail of the count distribution, which falsely predicts the presence of (higher-order) correlations. The related problems are lessened in experiments avoiding rapid changes of external conditions and long observation periods. In regard to the conflict between stationarity and efficiency, which demand short and long observation intervals, respectively, the parsimonious parametrization of the higher-order interactions appears as a particularly advantageous feature of our approach, since it helps to limit the length of the observation period required for reliable inferences.

Alternative parametrizations could also contribute to alleviate the consequences of instationarity. Thus one could treat the jump rate $\nu_+$ of the carrier process as a nuisance parameter and consider the normalized rates $\omega_n = \nu_n/\nu_+$ as the parameters of interest. Such would be appropriate particularly if only the general activation level $\nu_+$ is changing across time while the proportions between the rates $\nu_n$ remain unaffected. Stabilizing the variances of estimates may also be desirable, which would suggest the parametrization $\omega_n = (\nu_n/\nu_+)^{1/2}$.

For simplicity we assumed that the "coordinated" spikes that constitute a pattern in our model of population activity occur strictly simultaneously. Such precision of synchrony cannot be expected in reality. Intuitively, this would not hamper our theory to the degree to which the "jitter" is small enough compared to the width $h$ of the sampling bins, which was 5–20 milliseconds in most test simulations we looked at. In fact, this generates yet another conflict of interest between "small bins" which improve the performance of the estimators, and "large bins" which are less affected by jittered patterns (and, perhaps, by residual correlation across time resulting from spike after-effects).

What could be gained from analyses such as those proposed in this paper?

Positive statistical evidence for higher-order correlations in neuronal spike trains would be a novel result, potentially challenging conventional rate coding theories, as well as statements about the sufficiency of second-order correlations in capturing "essential aspects" of population activity in certain neuronal systems (e.g. Schneidman et al. [25], Shlens et al. [27]). However, spiking activity in a recurrent network will by necessity have a non-trivial statistical structure, potentially including also higher-order correlations. This is a consequence of network features, as for instance shared input and multiple feedback loops. Careful analysis of corresponding network models is necessary in order to correctly identify such structurally determined correlations, and not to mistake them for signatures of dynamic neuronal assemblies.

## Appendix A: Validity of the asymptotics

Our goal here is to investigate under what conditions the heuristic approximations of Section 3 can be expected to be valid. Of particular interest are consistency of estimates and the validity of the first-order (linear) approximation (3.2). We first discuss a directly tractable special case.

**A.1. A simple test case.** Let us consider the asymptotics of the estimate $\widehat{\nu}_+ = -h^{-1} \log \widehat{p}_0$ (cf. (4.8)) in some detail. The random variable $L\widehat{p}_0$ is binomially distributed with parameters $(L, p_0)$. Hence $E(\widehat{p}_0/p_0) = 1$, and $\mathrm{Var}(\widehat{p}_0/p_0) = (1-p_0)/(Lp_0) \to 0$ if $Lp_0 \to \infty$. It follows by a stochastic Taylor expansion of the logarithm that under this condition

$$\left(\frac{Lp_0}{1-p_0}\right)^{1/2} h\,(\widehat{\nu}_+ - \nu_+) = -\left(\frac{Lp_0}{1-p_0}\right)^{1/2} \log \frac{\widehat{p}_0}{p_0}$$
$$= -\left(\frac{Lp_0}{1-p_0}\right)^{1/2} \left(\frac{\widehat{p}_0}{p_0} - 1 + O_p\left(\frac{1-p_0}{Lp_0}\right)\right) = -B + o_p(1),$$

where $B = (p_0(1-p_0)/L)^{-1/2}(\widehat{p}_0 - p_0)$ is asymptotically standard normally distributed if and only if $Lp_0(1-p_0) \to \infty$. We conclude that $\widehat{\nu}_+$ *is asymptotically normally distributed, with mean $\nu_+$ and variance $h^{-2}(1-p_0)/(Lp_0)$, if and only if $Lp_0(1-p_0) \to \infty$.*

Let us yet explore the consequences of three natural assumptions:
(a) the intensity of the carrier ppp is bounded away from zero, $\nu_+^{-1} = O(1)$;
(b) the bin width $h$ stays bounded, $h = O(1)$;
(c) $\widehat{\nu}_+$ is consistent.
Consistency of $\widehat{\nu}_+$ means that its asymptotic variance tends to zero, $h^{-2}(1-p_0)/(Lp_0) = (e^{h\nu_+} - 1)/(Th) \to 0$. Since $\nu_+/T \leq (e^{h\nu_+} - 1)/(Th)$, this implies $\nu_+/T = o(1)$, which is thus a *necessary* condition for the consistency of $\widehat{\nu}_+$. Incidentally, the $h$-independence of the ultimate lower bound $\nu_+/T$ for the asymptotic variance shows that *choosing the bin width $h$ small is only partially helpful; given $\nu_+$ the rate of convergence can be effectively improved only by increasing the observation length $T$.* Another implication of consistency is $(Lp_0(1-p_0))^{-1} = o(1)$, i.e., *consistency of $\widehat{\nu}_+$ implies its asymptotic normality*. In fact, let us first suppose that $h\nu_+$ stays bounded away from zero. Then

$$(Lp_0(1-p_0))^{-1} = hT^{-1}e^{h\nu_+}/(1-e^{-h\nu_+}) = O(h^2 e^{h\nu_+}/(Th)),$$

and this tends to zero by (b) and (c). Similarly, if $h\nu_+ = o(1)$ then by (a) and (c)

$$hT^{-1}e^{h\nu_+}/(1-e^{-h\nu_+}) \sim h/(Th\nu_+) = \nu_+^{-2}(\nu_+/T) = o(1).$$

These considerations, along with the exponential dependence of the asymptotic variance on $h\nu_+ = -\log p_0$, point to the crucial role of the quantity $h\nu_+$ — which is evident, on the other hand, considering that $h\nu_+$ is the expected number of points of the carrier process in a bin. The present findings are corroborated in the following.

**A.2. General case.** Consistency of the estimates $\widehat{\nu}_n$ requires that the log-ecf, $\log \widehat{\gamma}_h$, is close in some sense to the logarithm of the common characteristic function of the bin counts $Z_l$, $\log \gamma_h$. In view of the expansions (3.1), (3.2), closeness largely depends on the order of magnitude of the random variables $h^{-1}\xi_h(\theta)$, which is immediate from (3.3):

$$\begin{aligned} E\,|h^{-1}\xi_h(\theta)|^2 &= (Th)^{-1}\left(|\gamma_h(\theta)|^{-2}-1\right) \\ &= (Th)^{-1}\left(\exp\left[2h\sum_k \nu_k(1-\cos(\theta k))\right]-1\right). \end{aligned}$$

The last expression is $\leq (Th)^{-1}(e^{4h\nu_+}-1)$ for every $\theta$, so on average (across $\theta$) we have by Jensen's inequality

$$\begin{aligned} \frac{e^{4h\nu_+}-1}{Th} &\geq \frac{1}{2\pi}\int_{-\pi}^{\pi} E\,|h^{-1}\xi_h(\theta)|^2\,d\theta &\text{(A.1)} \\ &\geq (Th)^{-1}\left(\exp\left[\frac{1}{2\pi}\int_{-\pi}^{\pi} 2h\sum_k \nu_k(1-\cos(\theta k))\,d\theta\right]-1\right) \\ &= \frac{e^{2h\nu_+}-1}{Th}. \end{aligned}$$

Let us remark here that the expression $\Xi_1 := (2\pi)^{-1}\int_{-\pi}^{\pi} E\,|h^{-1}\xi_h(\theta)|^2\,d\theta$ represents the quantity relevant to uniform consistency of all estimators of the form (2.5). This is (i) because to first order, one has for any linear functional $\lambda_c = \sum_k c_k \nu_k$ the approximation

$$\widehat{\lambda}_c - \lambda_c \approx \frac{1}{2\pi}\int_{-\pi}^{\pi} h^{-1}\xi_h(\theta)\,C(\theta)\,d\theta\,;$$

and (ii) because the last expression tends to 0 in mean square, uniformly across all $\lambda_c$ such that $||C||_2^2 = (2\pi)^{-1}\int_{-\pi}^{\pi}|C(\theta)|^2\,d\theta \leq 1$, if and only if $\Xi_1^{1/2} \to 0$. With $g(x) = (e^x-1)/x$, the lower bound in (A.1) becomes

$$\Xi_1 \geq g(2h\nu_+)\,2\nu_+/T \geq 2\nu_+/T.$$

Consequently, $\nu_+/T \to 0$ is a minimal requirement for the uniform consistency of the estimates; it is also sufficient if $h\nu_+$ stays bounded. Summarizing so far, *the following conditions are found to be necessary[7] for uniform consistency of the estimates:*

(C0)   $T \to \infty, \ h = O(1) \ (\text{so that } L = T/h \to \infty), \ \text{and } \nu_+/T \to 0.$

Going beyond consistency, let us now examine the (uniform) validity of the asymptotic (co-)variance approximations. These depend on the validity of the linear stochastic approximation (3.2), hence on the negligibility of the quadratic term (the first one omitted) in (3.1). *Its* expected mean square — which is relevant to the present question for the same reason as $\Xi_1$ was for uniform consistency — can be estimated as in the following lemma proven below.

**Lemma.** *Let $\Xi_2 = (2\pi)^{-1} \int_{-\pi}^{\pi} E \, |h^{-1}\xi_h(\theta)^2|^2 \, d\theta$. Then*

$$\Xi_2 \leq 3 \left( \frac{e^{4h\nu_+} - 1}{T} \right)^2 + 3 \frac{h}{T^3} \left( e^{8h\nu_+} - 1 \right) =: \epsilon_2 \,. \qquad (A.2)$$

Assuming (C0), under what additional conditions is $\epsilon_2$ of smaller order than the minimum order of the main term, $\Xi_1$? First, suppose that $h\nu_+$ stays bounded. Then $g(4h\nu_+)$ and $g(8h\nu_+)$ stay bounded as well, and the main term is of the order $\nu_+/T$. Writing $\epsilon_2$ in the form

$$\epsilon_2 = 3h^2 \left( \frac{4\nu_+}{T} g(4h\nu_+) \right)^2 + 3 \left( \frac{h}{T} \right)^2 \frac{8\nu_+}{T} g(8h\nu_+)$$

one finds that indeed $\epsilon_2 = o(\nu_+/T)$ under the assumed conditions. Thus the quadratic term is of smaller order than the main term, hence negligible.

Next, suppose that $h\nu_+$ becomes large. Then the main term $\Xi_1$ is at least of the order $e^{2h\nu_+}/(Th)$, by (A.1), while

$$\epsilon_2 \leq 3 \left( \frac{e^{2h\nu_+}}{Th} \right)^4 h^4 T^2 + 3 \frac{h}{T^3} \left( \frac{e^{2h\nu_+}}{Th} \right)^4 h^4 T^4 = 3 \left( \frac{e^{2h\nu_+}}{Th} \right)^4 \left( h^4 T^2 + h^5 T \right).$$

Hence, if $e^{2h\nu_+}/(Th) = o(1)$ as is necessary for uniform consistency in this case, then the condition $h^2 T = O(1)$ implies $\epsilon_2 = o\left(e^{2h\nu_+}/(Th)\right)$, and the quadratic term is negligible compared to the main term.

Putting everything together, we find that *the following conditions additionally to (C0) are decisive for the validity of the proposed asymptotics.*

(C1)   Either $h\nu_+ = O(1)$; or $h\nu_+ \to \infty, \ e^{2h\nu_+}/(Th) \to 0, \ \text{and } h^2 T = O(1).$

*Proof of lemma.* In order to prove (A.2), let us introduce the random variables $\zeta_l(\theta) = e^{i\theta Z_l}/\gamma(\theta)$. (Henceforth we supress $h$ whenever it appears as a parameter only.) Then

$$\xi(\theta)^2 = L^{-2} \sum\nolimits_{l_1, l_2} (\zeta_{l_1}(\theta) - 1)(\zeta_{l_2}(\theta) - 1),$$

---
[7] Boundedness of $h$ may not really be necessary, but is a natural condition in our setting.

and so because of $\overline{\zeta_l(\theta)} = \zeta_l(-\theta)$,

$$E|\xi(\theta)^2|^2 = L^{-4} \sum_{k_1, k_2, l_1, l_2} E(\zeta_{k_1}(\theta)-1)(\zeta_{k_2}(\theta)-1)(\zeta_{l_1}(-\theta)-1)(\zeta_{l_2}(-\theta)-1).$$

Since each single factor in the fourfold product has expectation zero and the $\zeta_l(\theta)$ are independent, only those terms contribute to the fourfold sum for which every two or all four indices are identical. Therefore,

$$E|\xi(\theta)^2|^2 = L^{-4} \sum_k E(\zeta_k(\theta)-1)^2 (\zeta_k(-\theta)-1)^2 \tag{A.3}$$
$$+ L^{-4} \sum_{k \neq l} E(\zeta_k(\theta)-1)^2 E(\zeta_l(-\theta)-1)^2$$
$$+ 2L^{-4} \sum_{k \neq l} E(\zeta_k(\theta)-1)(\zeta_k(-\theta)-1) E(\zeta_l(-\theta)-1)(\zeta_l(-\theta)-1)$$
$$= \frac{1}{L^3} E(\zeta(\theta)-1)^2 (\zeta(-\theta)-1)^2 + \frac{L-1}{L^3} E(\zeta(\theta)-1)^2 E(\zeta(-\theta)-1)^2$$
$$+ \frac{2(L-1)}{L^3} \left[E(\zeta(\theta)-1)(\zeta(-\theta)-1)\right]^2,$$

where $\zeta(\theta)$ denotes a generic random variable $\zeta_l(\theta)$.

Let us now evaluate the single expectations, in reverse order. Similarly as in (3.3) we get

$$E(\zeta(\theta)-1)(\zeta(-\theta)-1) = E\zeta(\theta)\zeta(-\theta) - 1 \tag{A.4}$$
$$= |\gamma(\theta)|^{-2} - 1 = \exp\left[2h \sum \nu_k (1-\cos(k\theta))\right] - 1 \leq e^{4h\nu_+} - 1.$$

Next, since

$$EX\,E\overline{X} = |EX|^2 \leq (E|X|)^2 \quad \text{and} \quad E(\zeta(\theta)-1)(\zeta(-\theta)-1) = E|\zeta(\theta)-1|^2,$$

$$E(\zeta(\theta)-1)^2 E(\zeta(-\theta)-1)^2 \leq \left(E|\zeta(\theta)-1|^2\right)^2 \tag{A.5}$$
$$= \left(\exp\left[2h \sum \nu_k (1-\cos(k\theta))\right] - 1\right)^2 \leq \left(e^{4h\nu_+} - 1\right)^2.$$

Last,

$$E(\zeta(\theta)-1)^2(\zeta(-\theta)-1)^2 \tag{A.6}$$
$$= E\left[\left(\frac{e^{i\theta Z}}{\gamma(\theta)} - 1\right)\left(\frac{e^{-i\theta Z}}{\gamma(-\theta)} - 1\right)\right]^2$$
$$= E\left[\frac{1}{|\gamma(\theta)|^2} + 1 - \frac{e^{i\theta Z}}{\gamma(\theta)} - \frac{e^{-i\theta Z}}{\gamma(-\theta)}\right]^2$$
$$= \left(\frac{1}{|\gamma(\theta)|^2} + 1\right)^2 - 4\left(\frac{1}{|\gamma(\theta)|^2} + 1\right) + \frac{\gamma(2\theta)}{\gamma(\theta)^2} + \frac{\gamma(-2\theta)}{\gamma(-\theta)^2} + \frac{2}{|\gamma(\theta)|^2}$$
$$= \frac{1}{|\gamma(\theta)|^4} - 1 + 2\left(\Re\frac{\gamma(2\theta)}{\gamma(\theta)^2} - 1\right).$$

But

$$\left| \Re \frac{\gamma(2\theta)}{\gamma(\theta)^2} \right| \leq \left| \frac{\gamma(2\theta)}{\gamma(\theta)^2} \right| = \exp\left[ h \sum \nu_k \left(\cos(2k\theta) - 1\right) - 2h \sum \nu_k \left(\cos(k\theta) - 1\right) \right]$$

$$= \exp\left[ h \sum \nu_k \left(\cos(2k\theta) - 2\cos(k\theta) + 1\right) \right]$$

$$\leq e^{h\nu_+/2},$$

the latter since $\cos(2k\theta) - 2\cos(k\theta) + 1 = 2\cos(k\theta)(1 - \cos(k\theta)) \leq 1/2$, so that in view of (A.6) we get

$$E\,(\zeta(\theta) - 1)^2(\zeta(-\theta) - 1)^2 \qquad (A.7)$$
$$\leq \exp\left[ 4h \sum \nu_k (1 - \cos(k\theta)) \right] - 1 + 2(e^{h\nu_+/2} - 1)$$
$$\leq 3(e^{8h\nu_+} - 1).$$

Putting everything together we finally obtain the bound

$$E\,|h^{-1}\xi(\theta)^2|^2 \leq \frac{3(L-1)}{h^2 L^3} \left(e^{4h\nu_+} - 1\right)^2 + \frac{3}{h^2 L^3}\left(e^{8h\nu_+} - 1\right)$$
$$\leq 3\left(\frac{e^{4h\nu_+} - 1}{T}\right)^2 + 3\frac{h}{T^3}\left(e^{8h\nu_+} - 1\right),$$

and the desired estimate (A.2) follows.

In passing we note the related lower estimate

$$\frac{1}{2\pi} \int_{-\pi}^{\pi} E\,|h^{-1}\xi(\theta)^2|^2\,d\theta \geq 2(1 - h/T)\left(\frac{e^{2h\nu_+} - 1}{T}\right)^2. \qquad (A.8)$$

It can be established similarly as earlier on noting first that by (A.3) and (A.4)

$$E\,|\xi(\theta)^2|^2 \geq \frac{2(L-1)}{L^3}\left[E\,(\zeta(\theta) - 1)(\zeta(-\theta) - 1)\right]^2$$
$$= \frac{2(L-1)}{L^3}\left(\exp\left[2h\sum \nu_k(1 - \cos(k\theta))\right] - 1\right)^2,$$

then using the convexity of the function $f(x) = (e^x - 1)^2$ $(x \geq 0)$ along with Jensen's inequality.

### Appendix B: Evaluation of asymptotic (co-)variances

All asymptotic (co-)variances of interest to us are of the form $T^{-1}J(n_1, n_2, \alpha_1, \alpha_2)$ for some $n_j \geq 0$ and $\alpha_j \in \{0, 1\}$ $(j = 1, 2)$, where

$$J(n_1, n_2, \alpha_1, \alpha_2) = \frac{1}{(2\pi)^2}\int_{-\pi}^{\pi}\int_{-\pi}^{\pi}\frac{\Gamma_h(\theta_1, \theta_2)}{(e^{i\theta_1} - 1)^{\alpha_1}(e^{i\theta_2} - 1)^{\alpha_2}}\frac{d\theta_1}{e^{i\theta_1 n_1}}\frac{d\theta_2}{e^{i\theta_2 n_2}}. \qquad (B.1)$$

This follows from (3.6) and the subsequent remark on noting that for an individual rate estimate $\widehat{\nu}_n$ we have $C(\theta) = e^{-i\theta n}$, while for the tail estimate $\widehat{\rho}_m$ we have $C(\theta) = e^{-i\theta m}/(1 - e^{-i\theta}) = e^{-i\theta(m-1)}/(e^{i\theta} - 1)$. In particular,

$$\operatorname{ascov}(\widehat{\nu}_{n_1}, \widehat{\nu}_{n_2}) = T^{-1} J(n_1, n_2, 0, 0) = T^{-1} \Omega_{n_1, n_2}, \tag{B.2}$$

$$\operatorname{ascov}(\widehat{\rho}_{m_1}, \widehat{\rho}_{m_2}) = T^{-1} J(m_1-1, m_2-1, 1, 1) = T^{-1} \Sigma_{m_1, m_2}, \tag{B.3}$$

$$\operatorname{ascov}(\widehat{\rho}_m, \widehat{\nu}_n) = T^{-1} J(m-1, n, 1, 0). \tag{B.4}$$

To evaluate the double integral (B.1) we utilize the explicit representation (5.1) of the kernel $\Gamma_h(\theta_1, \theta_2)$, from which it follows that the integrand of (B.1) depends on the arguments $\theta_j$ only via the variables $e^{i\theta_j}$ $(j = 1, 2)$. The substitutions $\theta_j \mapsto e^{i\theta_j} = z_j$ then give

$$J(n_1, n_2, \alpha_1, \alpha_2) = \frac{1}{(2\pi i)^2} \int_{|z_1|=1} \int_{|z_2|=1} \psi_{\alpha_1, \alpha_2}(z_1, z_2) \frac{dz_1}{z_1^{n_1+1}} \frac{dz_2}{z_2^{n_2+1}} \tag{B.5}$$

where

$$\psi_{\alpha_1, \alpha_2}(z_1, z_2) = \frac{\exp\left[h \sum_{n=1}^{N} \nu_n (z_1^n - 1)(z_2^n - 1)\right] - 1}{h (z_1 - 1)^{\alpha_1} (z_2 - 1)^{\alpha_2}}. \tag{B.6}$$

No matter what $\alpha_j \in \{0, 1\}$, the function $\psi_{\alpha_1, \alpha_2}(z_1, z_2)$ is a holomorphic (in fact, entire) function of the two complex variables $z_1, z_2$. Therefore, the bivariate version of Cauchy's formula applies, and we obtain the following evaluation of (B.1) (and hence, of the asymptotic covariances),

$$J(n_1, n_2, \alpha_1, \alpha_2) = \frac{1}{n_1! n_2!} \left. \frac{\partial^{n_1+n_2} \psi_{\alpha_1, \alpha_2}(z_1, z_2)}{\partial z_1^{n_1} \partial z_2^{n_2}} \right|_{z_1=z_2=0}. \tag{B.7}$$

Let us give a few examples. In the test case of Section A the asymptotic variance of $\widehat{\nu}_+ = \widehat{\rho}_1$ times $T$ was found to be $(e^{h\nu_+} - 1)/h$. The same expression should obtain for $J(0, 0, 1, 1) = \psi_{1,1}(0, 0)$, and in fact, it does. For $T \operatorname{asvar}(\widehat{\nu}_1)$ one finds $J(1, 1, 0, 0) = e^{h\nu_+}(\nu_1 + h\nu_1^2)$, and similar formulae result for $\widehat{\nu}_n$, $n \geq 2$. Examples of covariances are

$$\Omega_{1,2} = T \operatorname{ascov}(\widehat{\nu}_1, \widehat{\nu}_2) = e^{h\nu_+} h \nu_1 \left(\nu_2 - \nu_1 - \frac{1}{2} h \nu_1^2\right),$$

$$T \operatorname{ascov}(\widehat{\rho}_1, \widehat{\nu}_1) = e^{h\nu_+} \nu_1.$$

Higher-order formulae are a mess, but can be obtained using a computer algebra system.

**Acknowledgment**

We thank the referee for his/her insightful comments and constructive suggestions, and for pointing out an error in the first version of the manuscript.

**Added Note**

The following two articles relevant to the contents of the present paper came to the attention of the authors immediately before publication:

Buchmann, B. and Grübel, R. (2003). Decompounding: an estimation problem for Poisson random sums. *Ann. Statist.* **31**, 1054-1074.

van Es, B., Gugushvili, S. and Spreij, P. (2007). A kernel type nonparametric density estimator for decompounding. *Bernoulli* **13**, 672-694.

See also the references cited in these works.

**References**

[1] Abeles, M. (1991). *Corticonics: Neural circuits of the cerebral cortex.* Cambridge University Press.
[2] Ahlfors, L. V. (1979). *Complex analysis, 3rd ed.* McGraw-Hill, New York. MR0510197
[3] Amarasingham, A., Chen, T-L., Harrison, M. T. and Sheinberg, D. L. (2006). Spike count reliability and the Poisson hypothesis. *J. Neurosci.* **26**, 801-809.
[4] Baker, S. N. and Lemon, R. N. (2000). Precise spatiotemporal repeating patterns in monkey primary and supplementary motor areas occur at chance levels. *J. Neurophysiology* **84**, 1770-1780.
[5] Borodin, A., Okounkov, A. and Olshanski, G. (2000). Asymptotics of Plancherel measures for symmetric groups. *J. Amer. Math. Soc.* **13**, 481–515. MR1758751
[6] Brown, E. N., Kass, R. E. and Mitra, P. P. (2004). Multiple neural spike train data analysis: State-of-the-art and future challenges. *Nature Neurosci.* **7**, 456–461.
[7] Daley, D.J. and Vere-Jones, D. (2002). *An introduction to the theory of point processes: Elementary theory and methods, Vol. 1.* Springer, New York. MR0950166
[8] Diesmann, M., Gewaltig, M.-O. and Aertsen, A. (1999). Stable propagation of synchronous spiking in cortical neural networks. *Nature* **402**, 529-533.
[9] Grün, S., Diesmann, M. and Aertsen, A. (2002a). Unitary events in multiple single-neuron spiking activity: I. Detection and significance. *Neural Computation* **14**, 43-80.
[10] Grün, S., Diesmann, M. and Aertsen, A. (2002b). Unitary events in multiple single-neuron spiking activity: II. Nonstationary data. *Neural Computation* **14**, 81-119.
[11] Gütig, R., Aertsen, A. and Rotter, S. (2002). Statistical significance of coincident spikes: Count-based versus rate-based statistics. *Neural Computation* **14**, 121-153.
[12] Gütig R., Aertsen, A. and Rotter, S. (2003). Analysis of higher-order neuronal interactions based on conditional inference. *Biol. Cybern.* **88**, 352-359.
[13] Hebb, D. O. (1949). *The organization of behavior: A neuropsychological theory.* Wiley, New York.

[14] Jacquet, P. and Szpankowski, W. (1998). Analytical de-Poissonization and its applications. *Theoret. Comput. Sci.* **201**, 1–62. MR1625392

[15] Jongbloed, G., van der Meulen, F. H. and van der Vaart, A. W. (2005). Nonparamatric inference for Lévy driven Ornstein-Uhlenbeck processes. *Bernoulli* **11**, 759-791. MR2172840

[16] Kuhn, A., Rotter, S. and Aertsen, A. (2002). Correlated input spike trains and their effects on the response of the leaky integrate-and-fire neuron. *Neurocomputing* **44-46**, 121-126.

[17] Kuhn, A., Aertsen, A. and Rotter, S. (2003). Higher-order statistics of input ensembles and the response of simple model neurons. *Neural Computation* **15**, 67-101.

[18] Lee, A. K., Manns, I. D., Sakmann, B. and Brecht, M. (2006). Whole-cell recordings in freely moving rats. *Neuron* **51**, 399–407.

[19] von der Malsburg, C. (1995). Binding in models of perception and brain function. *Current Opinion in Neurobiology* **5**, 520-526.

[20] Martignon, L., von Hasseln, H., Grün, S., Aertsen, A. and Palm, G. (1995). Detecting higher-order interactions among the spiking events in a group of neurons. *Biol. Cybern.* **73**, 69–81.

[21] Martignon, L., Deco, G., Laskey, K., Diamond, M., Freiwald, W. and Vaadia, E. (2000). Neural coding: higher-order temporal patterns in the neurostatistics of cell assemblies. *Neural Computation* **12**, 2621–2653.

[22] Nakahara, H. and Amari, S. (2002). Information-geometric measure for neural spikes. *Neural Computation* **14**, 2269–2316.

[23] Perkel, D. H., Gerstein, G. L., Smith, M. S. and Tatton, W. G. (1975). Nerve-impulse patterns: a quantitative display technique for three neurons. *Brain Research* **100**, 271-296.

[24] Prakasa Rao, B. L. S. (1987). *Asymptotic theory of statistical inference.* Wiley, New York. MR0874342

[25] Schneidman, E., Berry, M. J. 2nd, Segev, R. and Bialek, W. (2006). Weak pairwise correlations imply strongly correlated network states in a neural population. *Nature* **440**, 1007–1012.

[26] Shadlen, M. N. and Newsome, W. T. (1998). The variable discharge of cortical neurons: implications for connectivity, computation, and information coding. *J. Neurosci.* **18**, 3870–3896.

[27] Shlens, J., Field, G. D., Gauthier, J. L., Grivich, M. I., Petrusca, D., Sher, A., Litke, A. M. and Chichilnisky, E.J. (2006). The structure of multi-neuron firing patterns in primate retina. *J. Neurosci.* **26**, 8254–8266.